\documentclass{article}
\usepackage{amssymb,latexsym,amsmath,amsthm,mathrsfs}
\usepackage{graphicx}
\newcommand{\comment}[1]{}

\begin{document}
\title{Expansion of the infinite product
$(1-x)(1-xx)(1-x^3)(1-x^4)(1-x^5)(1-x^6)\; \textrm{etc.}$
into a simple series\footnote{Presented to the St. Petersburg Academy
on August 14, 1775. Originally published as
{\em Evolutio producti infiniti
$(1-x)(1-xx)(1-x^3)(1-x^4)(1-x^5)(1-x^6)\; \textrm{etc.}$
in seriem simplicem},
Acta academiae scientiarum Petropolitanae \textbf{1780} (1783), part I, 
47--55.
E541 in the Enestr{\"o}m index.
Translated from the Latin by Jordan Bell,
Department of Mathematics, University of Toronto, Toronto, Ontario, Canada.
Email: jordan.bell@gmail.com}}
\author{Leonhard Euler}
\date{}
\maketitle

1. Putting
\[
s=(1-x)(1-xx)(1-x^3)(1-x^4)\;\textrm{etc.}
\]
it is easy to see that
\[
s=1-x-xx(1-x)-x^3(1-x)(1-xx)-x^4(1-x)(1-xx)(1-x^3)-\textrm{etc.};
\]
now, as this series is infinite, if all its terms are expanded one searches
for what kind of series will be produced according to simple powers
of $x$. Then since the first two terms $1-x$ have already been expanded,
one may write the letter $A$ in place of all the others, so that it will be
\[
s=1-x-A
\]
and hence
\[
A=xx(1-x)+x^3(1-x)(1-xx)+x^4(1-x)(1-xx)(1-x^3)+\textrm{etc.}
\]

2. Since all these terms have the common factor $1-x$, with that expanded
all the terms will be separated into two parts, which we will represent
thus
{\footnotesize
\[
\begin{array}{rlllll}
A=&xx&+x^3(1-xx)&+x^4(1-xx)(1-x^3)&+x^5(1-xx)(1-x^3)(1-x^4)&+\textrm{etc.}\\
&-x^3&-x^4(1-xx)&-x^5(1-xx)(1-x^3)&-x^6(1-xx)(1-x^3)(1-x^4)&-\textrm{etc.}
\end{array}
\]
}
Now here, the two parts having the same power of $x$ are added together, and the following form will result for $A$
\[
A=xx-x^5-x^7(1-xx)-x^9(1-xx)(1-x^3)-x^{11}(1-x^2)(1-x^3)(1-x^4)-\textrm{etc.}
\]
where now the first two terms $xx-x^5$ have been expanded;
and the following proceed in the powers $x^7,x^9,x^{11},x^{13},x^{15}$, whose
exponents increase by two.

3. Now let us put in a similar way as before
\[
A=xx-x^5-B
\]
so that it will be
\[
B=x^7(1-xx)+x^9(1-xx)(1-x^3)+x^{11}(1-xx)(1-x^3)(1-x^4)+\textrm{etc.},
\]
all of whose terms have the common factor $1-xx$. With this expanded, all the
terms are separated into two parts, so that it follows that
{\footnotesize
\[
\begin{array}{rlllll}
B=&x^7&+x^9(1-x^3)&+x^{11}(1-x^3)(1-x^4)&+x^{13}(1-x^3)(1-x^4)(1-x^5)&+\textrm{etc.}\\
&-x^9&-x^{11}(1-x^3)&-x^{13}(1-x^3)(1-x^4)&-x^{15}(1-x^3)(1-x^4)(1-x^5)&-\textrm{etc.}
\end{array}
\]
}
Here again the two terms which have the same power of $x$ in front are
collected together, and this will produce
\[
B=x^7-x^{12}-x^{15}(1-x^3)-x^{18}(1-x^3)(1-x^4)-x^{21}(1-x^3)(1-x^4)(1-x^5)-\textrm{etc.}
\]
where now all the powers of $x$ increase by three.

4. Now it is in turn put
\[
B=x^7-x^{12}-C,
\]
so that it will be
\[
C=x^{15}(1-x^3)+x^{18}(1-x^3)(1-x^4)+x^{21}(1-x^3)(1-x^4)(1-x^5)+\textrm{etc.},
\]
and by expanding the factor $1-x^3$, all the terms are resolved into two parts
and it will become
{\footnotesize
\[
\begin{array}{rlllll}
C=&x^{15}&+x^{18}(1-x^4)&+x^{21}(1-x^4)(1-x^5)&+x^{24}(1-x^4)(1-x^5)(1-x^6)&+\textrm{etc.}\\
&-x^{18}&-x^{21}(1-x^4)&-x^{24}(1-x^4)(1-x^5)&-x^{27}(1-x^4)(1-x^5)(1-x^6)&-\textrm{etc.},
\end{array}
\]}
where once again the members with the same power of $x$ in front joined
together yield
\[
C=x^{15}-x^{22}-x^{26}(1-x^4)-x^{30}(1-x^4)(1-x^5)-x^{34}(1-x^4)(1-x^5)(1-x^6)-\textrm{etc.},
\]
where the powers in front increase by four.

5. Let it be put
\[
C=x^{15}-x^{22}-D,
\]
so that it will be
\[
D=x^{26}(1-x^4)+x^{30}(1-x^4)(1-x^5)+x^{34}(1-x^4)(1-x^5)(1-x^6)+\textrm{etc.}
\]
By expanding the factor $1-x^4$, the terms will
be separated into two thus
{\footnotesize
\[
\begin{array}{rlllll}
D=&x^{26}&+x^{30}(1-x^5)&+x^{34}(1-x^5)(1-x^6)&+x^{38}(1-x^5)(1-x^6)(1-x^7)&+\textrm{etc.}\\
&-x^{30}&-x^{34}(1-x^5)&-x^{38}(1-x^5)(1-x^6)&-x^{42}(1-x^5)(1-x^6)(1-x^7)&-\textrm{etc.}
\end{array}
\]
}
Now with the two terms joined together as before, it will turn out as
\[
D=x^{26}-x^{35}-x^{40}(1-x^5)-x^{45}(1-x^5)(1-x^6)-x^{50}(1-x^5)(1-x^6)(1-x^7)-\textrm{etc.}
\]
So here the powers of $x$ increase by five.

6. Let it be put
\[
D=x^{26}-x^{35}-E,
\]
so that it will be
\[
E=x^{40}(1-x^5)+x^{45}(1-x^5)(1-x^6)+x^{50}(1-x^5)(1-x^6)(1-x^7)+\textrm{etc.},
\]
and by doing the resolution into two parts as before it will yield
{\footnotesize
\[
\begin{array}{rlllll}
E=&x^{40}&+x^{45}(1-x^6)&+x^{50}(1-x^6)(1-x^7)&+x^{55}(1-x^6)(1-x^7)(1-x^8)&+\textrm{etc.}\\
&-x^{45}&-x^{50}(1-x^6)&-x^{55}(1-x^6)(1-x^7)&-x^{60}(1-x^6)(1-x^7)(1-x^8)&-\textrm{etc.}
\end{array}
\]
}
Indeed, with the two terms joined together this will produce
\[
E=x^{40}-x^{51}-x^{57}(1-x^6)-x^{63}(1-x^6)(1-x^7)-x^{69}(1-x^6)(1-x^7)(1-x^8)-\textrm{etc.},
\]
where all the powers of $x$ increase by six.

7. The rule by which these operations are to be continued is clear enough
if the first terms of each of the letters $A,B,C,D,E$ etc. are substituted
in order, and we will find the following form for the sought series
\[
s=1-x,-xx+x^5,+x^7-x^{12},-x^{15}+x^{22},+x^{26}-x^{35},-x^{40}+x^{51},+\textrm{etc.}
\]
Here the question reduces to the following: that the order be defined
by which the exponents of the powers of $x$ can be extended continually further,
since it is now clear enough from the operations that have been done that the
signs $+$ and $-$ of the terms follow each other alternately in pairs.

8. Therefore, to inquire into this law, let us see how these numbers arise
in all the letters. To this end let us arrange the first
terms in the first
forms given for each of
the letters:
\[
\begin{array}{ll|llll}
A&=xx(1-x)&7&=3+4&=3+1+3&=3+1+1+2\\
B&=x^7(1-xx)&15&=4+11&=4+2+9&=4+2+2+7\\
C&=x^{15}(1-x^3)&26&=5+21&=5+3+18&=5+3+3+15\\
D&=x^{26}(1-x^4)&40&=6+34&=6+4+30&=6+4+4+26\\
E&=x^{40}(1-x^5)&57&=7+50&=7+5+45&=7+5+5+40\\
&\textrm{etc.}&&&\textrm{etc.}&
\end{array}
\]
Here for instance, from the expansion of the letter $A$ we have seen
that the number $7$
arises from the sum $3+4$, and indeed $4$ arises from $1+3$
and too $3$ from $1+2$,\footnote{Translator: I don't understand the right hand
column in this table. I think it is just a way of explaining the second order
arithmetic progression, and does not refer to the way the numbers $7,15,26$, etc. actually
arise in Euler's proof. For instance, in \S 2, $-x^7(1-xx)$ comes from
adding the terms $-x^4(1-xx)$ and $x^4(1-xx)(1-x^3)$, thus $7=3+4$. And
$-x^4(1-xx)$ is one of the terms we get when we expand the $1-x$ in
$+x^3(1-x)(1-xx)$, thus $4=1+3$. But I don't see what
$3=1+2$ refers to. Another example is 40: I see where $40=6+34$ comes from,
but I don't see why 34 is being written as $4+30$ instead of $7+27$, which
is how 34 arises in \S 4.} which therefore gives the resolution
\[
7=3+4=3+1+3=3+1+1+2.
\]
And the same order is seen in the following letters, where the last numbers
proceed in the order $2,7,15,26,40$.

9. From this it is clear that the differences of the numbers $2,7,15,26,40,
57$ etc. constitute an arithmetic progression, whence the general term
of these numbers will be\footnote{Translator: Let
$\Delta a_n=a_{n+1}-a_n$ and $\Delta^{k+1}a_n=\Delta^k a_{n+1}-\Delta^k a_n$.
Then
\[
a_n=\sum_{k=0}^\infty \binom{n}{k} \Delta^k a_0.
\]
For $a_0=2,a_1=7,a_2=15$, Euler gives here $a_{n-1}$.}
\[
2+5(n-1)+\frac{3(n-1)(n-2)}{1\cdot 2}=\frac{3nn+n}{2}.
\]
Moreover, the exponents which precede these were $1,5,12,22,35,51$ etc. and
the differences with these are $1,2,3,4,5,6$
etc. and in general the number $n$ itself, so that the exponent which
precedes the formula $\frac{3nn+n}{2}$ will be 
\[
\frac{3nn-n}{2}.
\]

10. Now we understand perfectly the simple series that is found equal to
the given infinite product
\[
(1-x)(1-xx)(1-x^3)(1-x^4) \; \textrm{etc.},
\] 
For since the series has been found to be
\[
s=1-x^1-x^2+x^5+x^7-x^{12}-x^{15}+x^{22}+x^{26}-x^{35}-x^{40}+x^{51}
+\textrm{etc.},
\]
we are now certain that in it no other powers of occur occur unless their
exponents are contained in the general formula $\frac{3nn \pm n}{2}$,
and indeed even more precisely, if $n$ is an odd number then the two terms
arising from it will have the sign $-$ and those which arise from even
numbers the sign $+$.

\begin{center}
{\Large Another investigation of the same series}
\end{center}

11. The same series proceeding according to powers of $x$ can also be
investigated
in the following way. For as
\[
s=1-x-xx(1-x)-x^3(1-x)(1-xx)-x^4(1-x)(1-x^2)(1-x^3)-\textrm{etc.},
\]
expand the second term $-xx(1-x)$, so it becomes
\[
s=1-x-xx+x^3-x^3(1-x)(1-xx)-x^4(1-x)(1-xx)(1-x^3)-\textrm{etc.},
\]
and let
\[
s=1-x-xx+A,
\]
so that it is
\[
A=x^3-x^3(1-x)(1-xx)-x^4(1-x)(1-xx)(1-x^3)-\textrm{etc.}
\]
Now separate each term into two parts by expanding the factor $1-x$, 
which will produce
{\footnotesize
\[
\begin{array}{llllll}
A=&x^3&-x^3(1-xx)&-x^4(1-xx)(1-x^3)&-x^5(1-x^2)(1-x^3)(1-x^4)&-\textrm{etc.}\\
&&+x^4(1-xx)&+x^5(1-xx)(1-x^3)&+x^6(1-x^2)(1-x^3)(1-x^4)&+\textrm{etc.}
\end{array}
\]}
Here again, adding together all the pairs of terms having the same power of $x$
will yield
\[
A=x^5+x^7(1-xx)+x^9(1-xx)(1-x^3)+x^{11}(1-x^2)(1-x^3)(1-x^4)+\textrm{etc.}
\]

12. Again here, expand the second term, which will give
\[
A=x^5+x^7-x^9+x^9(1-xx)(1-x^3)+x^{11}(1-x^2)(1-x^3)(1-x^4)+\textrm{etc.}
\]
Now put
\[
A=x^5+x^7-B,
\]
so that
\[
B=x^9-x^9(1-xx)(1-x^3)-x^{11}(1-xx)(1-x^3)(1-x^4)-\textrm{etc.};
\]
if each factor $1-xx$ in this is expanded, we will obtain
{\footnotesize
\[
\begin{array}{llllll}
B=&x^9&-x^9(1-x^3)&-x^{11}(1-x^3)(1-x^4)&-x^{13}(1-x^3)(1-x^4)(1-x^5)&-\textrm{etc.}\\
&&+x^{11}(1-x^3)&+x^{13}(1-x^3)(1-x^4)&+x^{15}(1-x^3)(1-x^4)(1-x^5)&+\textrm{etc.},
\end{array}
\]}
then indeed by combining the 
two terms we get
\[
B=x^{12}+x^{15}(1-x^3)+x^{18}(1-x^3)(1-x^4)+x^{21}(1-x^3)(1-x^4)(1-x^5)+\textrm{etc.}
\]

13. Expand likewise the second term and put
\[
B=x^{12}+x^{15}-C
\]
and it will be
\[
C=x^{18}-x^{18}(1-x^3)(1-x^4)-x^{21}(1-x^3)(1-x^4)(1-x^5)-\textrm{etc.}
\]
Now expand the terms according to the factor $1-x^3$, and this will make
{\footnotesize
\[
\begin{array}{llllll}
C=&x^{18}&-x^{18}(1-x^4)&-x^{21}(1-x^4)(1-x^5)&-x^{24}(1-x^4)(1-x^5)(1-x^6)&-\textrm{etc.}\\
&&+x^{21}(1-x^4)&+x^{24}(1-x^4)(1-x^5)&+x^{27}(1-x^4)(1-x^5)(1-x^6)&+\textrm{etc.}
\end{array}
\]}
Then combining the two terms will make
\[
C=x^{22}+x^{26}(1-x^4)+x^{30}(1-x^4)(1-x^5)+x^{34}(1-x^4)(1-x^5)(1-x^6)+\textrm{etc.}
\]

14. Here, now expanding again the second term, put
\[
C=x^{22}+x^{26}-D
\]
and it will be
\[
D=x^{30}-x^{30}(1-x^4)(1-x^5)-x^{34}(1-x^4)(1-x^5)(1-x^6)-\textrm{etc.},
\]
where expanding the factor $1-x^4$ will yield
{\footnotesize
\[
\begin{array}{llllll}
D=&x^{30}&-x^{30}(1-x^5)&-x^{34}(1-x^5)(1-x^6)&-x^{38}(1-x^5)(1-x^6)(1-x^7)&-\textrm{etc.}\\
&&+x^{34}(1-x^5)&+x^{38}(1-x^5)(1-x^6)&+x^{42}(1-x^5)(1-x^6)(1-x^7)&+\textrm{etc.}
\end{array}
\]}
Combining the two terms gives
\[
D=x^{35}+x^{40}(1-x^5)+x^{45}(1-x^5)(1-x^6)+x^{50}(1-x^5)(1-x^6)(1-x^7)+\textrm{etc.}
\]

15. Next, expand the second term and put
\[
D=x^{35}+x^{40}-E,
\]
and it will be
\[
E=x^{45}-x^{45}(1-x^5)(1-x^6)-x^{50}(1-x^5)(1-x^6)(1-x^7)-\textrm{etc.}
\]
and expanding according to the factor $1-x^5$ 
will make
{\footnotesize
\[
\begin{array}{llllll}
E=&x^{45}&-x^{45}(1-x^6)&-x^{50}(1-x^6)(1-x^7)&-x^{55}(1-x^6)(1-x^7)(1-x^8)&-\textrm{etc.}\\
&&+x^{50}(1-x^6)&+x^{55}(1-x^6)(1-x^7)&+x^{60}(1-x^6)(1-x^7)(1-x^8)&+\textrm{etc.}
\end{array}
\]}
and collecting the two terms together elicits
\[
E=x^{51}+x^{57}(1-x^6)+x^{63}(1-x^6)(1-x^7)+x^{69}(1-x^6)(1-x^7)(1-x^8)+\textrm{etc.}
\]

16. Therefore, having found the values of the letters $A,B,C,D,E$ etc.,
if they are all successively substituted, this series will result
\[
1-x-xx,+x^5+x^7,-x^{12}-x^{15},+x^{22}+x^{26},-x^{35}-x^{40},+\textrm{etc.}
\]
It is easy to see the order of the exponents here. For since 
in the original expressions for the letters $A,B,C,D,E$ etc. the first
simple terms were $x^3,x^9,x^{18},x^{30},x^{45}$ etc., whose exponents are
clearly three times the triangular numbers; generally for the number $n$ the
exponent will be $\frac{3nn+3n}{2}$. 
These terms follow two powers of $x$ preceding by the same difference $n$,
whence by subtracting the number $n$ twice from this formula, the two
powers occurring in the sought series will be formed, whose exponents
will be
\[
\frac{3nn+n}{2} \quad \textrm{and} \quad \frac{3nn-n}{2}.
\]

17. Then on the other hand, it is apparent that the series
\[
s=1-x-xx+x^5+x^7-x^{12}-x^{15}+x^{22}+x^{26}-x^{35}-\textrm{etc.}
\]
continued to infinity has infinitely many factors, which will namely be
\[
1-x,\quad 1-xx,\quad 1-x^3,\quad 1-x^4,\quad 1-x^5\quad \textrm{etc.}
\]
so that, if it is first divided by $1-x$, then the quotient by $1-xx$,
and this quotient next by $1-x^3$, and this division were carried on thus
to infinity, the ultimate quotient resulting should be equal to unity.

18. So if this equation extending to infinity were given
\[
1-x-xx+x^5+x^7-x^{12}-x^{15}+x^{22}+x^{26}-x^{35}-\textrm{etc.}=0,
\]
all of its roots could easily be assigned. For the first root will be $x=1$,
then the two square roots of unity, then indeed the three cube roots of unity,
next the four fourth roots of unity, and similarly the five fifth roots
of unity, and so on, among which then unity will occur infinitely many times;
and indeed $-1$ will appear when the root of an even power is to be extracted.

\end{document}